\documentclass{amsart}\usepackage{amsthm,amssymb,amsmath,dsfont,bigdelim}
\usepackage{graphicx}
\usepackage{float}
\usepackage[arc,all]{xy}
\usepackage{longtable}
\newtheorem{thm}{Theorem}[section]
\newtheorem{cor}[thm]{Corollary}
\newtheorem{lem}[thm]{Lemma}
\newtheorem{prop}[thm]{Proposition}
\theoremstyle{definition}
\newtheorem{defn}[thm]{Definition}
\newtheorem{exa}[thm]{Example}

\newtheorem*{proof M}{\textbf{Proof of  Main Theorem}}
\theoremstyle{remark}

\numberwithin{equation}{section}

\begin{document}
\title[commutativity degree]{On the commutativity degree of a finite-dimensional  Lie algebra}%
\author{Afsaneh Shamsaki}%
\address{Department of Mathematics, Faculty of Mathematical Sciences,
Ferdowsi University of Mashhad, Mashhad, Iran }
\email{Shamsaki.afsaneh@yahoo.com}%
\author{Ahmad Erfanian }
\email{erfanian@um.ac.ir }
\address{Department of Mathematics, Faculty of Mathematical Sciences,
Ferdowsi University of Mashhad, Mashhad, Iran }
\author{Mohsen Parvizi}
\email{parvizi@um.ac.ir }
\address{Department of Mathematics, Faculty of Mathematical Sciences,
Ferdowsi University of Mashhad, Mashhad, Iran }
\begin{abstract}
In this paper, we introduce the commutativity degree of a finite-dimensional Lie algebra over a finite field and determine  upper and lower bounds for it. Moreover, we study some relations between the notion of commutativity degree and known concepts in Lie algebras.
\end{abstract}
\maketitle
\section{introduction}
In the past decades, the relation between probability theory and other branches of mathematics has attracted some authors. For instance, the relation between group theory and probability theory is considered and  the commutativity degree of a finite group $ G $ denoted by $d(G)$ defined as the probability
that two randomly chosen elements of $ G $ commute. More precisely 
\begin{equation*}
d(G)=\dfrac{| \lbrace (x, y)\in G\times G \mid xy=yx \rbrace |}{| G|^{2}}.
\end{equation*}
Gustafson  \cite{Gus} was the first one who obtained the upper bound $\frac{5}{8}$ for it. 
Later, Lescot \cite{Les} showed that if $ G $ is a group such that $ d(G)> \frac{1}{2}, $ then $ G $ is nilpotent. He studied the relation between the commutativity degree and the concept of isoclinism and proved that two isoclinic finite groups have the same commutativity degree, but the converse is not true in general. One of his results is the computation of  the commutativity degree of  $ D_{2n} $ and $ Q_{2^{n+1}} $.\\
Recently, Nath \cite{Nath} found the commutativity degree of finite $ p $-groups with derived subgroups of order $ p $ such that $ G^{\prime}\subseteq Z(G) $ and also proved that for every positive integer $ k $ there exists a family of groups with the  commutativity degree of $ \frac{1}{k} $. 
In \cite{Ring}, by Haar measure the notion of commutativity degree was extended for infinite compact groups. Moreover, the commutativity degree of Moufang loops was defined and investigated by Ahmadidelir \cite{Ahmadi}.
\\
In spite of the fact that two concepts ``Lie algebra'' and ``group'' are totally different and they have different structures, but sometimes we may observe similar behaviours. In this paper, we are going to define the commutativity degree of a finite-dimensional Lie algebra over the field $ \mathbb{F}_{q} $ by following the same lines in defining it in group theory. As this notion originates in group theory, one may expect a considerable similarity between the notion of commutativity degree in groups and Lie algebras, but there are some important structural differences either in the definitions or in techniques. For instance, the role of conjugacy classes in computing the commutativity degree of a group is very important, but there is no such definition in Lie algebras.   
\section{Preliminaries}
In the following, we recall the definition of a Lie algebra and   state some  necessary concepts and  terminologies which are used in the rest of the paper.\\
Let $ \mathbb{F} $ be  a field. A Lie algebra over the field $ \mathbb{F} $  is an $ \mathbb{F} $-vector space $ L, $ together with a bilinear map called the Lie bracket $ L\times L\rightarrow L $ with $ (x, y)\mapsto [x, y], $ satisfying the following properties:
\begin{align*}
&[x, x]=0\quad \text{for all} \quad x\in L, \quad & (L1)\cr
&[x, [y, z]]+[z,[x, y]]+[y,[z, x]]=0 \quad \text{for all} \quad x, y, z\in L. \quad &(L2)
\end{align*} 
The Lie bracket $ [x, y] $ is often referred to as  the commutator of $ x $ and $ y. $ The condition $ L2 $ is known as the Jacobi identity. Since the Lie bracket $ [-, -] $ is bilinear, we have $ 0=[x+y, x+y]=[x, x]+[x, y]+[y, x]+[y, y]. $ Hence the condition $ L1 $ implies $ [x, y]=-[y, x]$ for all $ x, y\in L. $  As an example of Lie algebras, if $ [x, y]=0 $ for all $ x, y \in L, $ then $ L $ is called an abelian Lie algebra.  In this paper, abelian $ n $-dimensional Lie algebras are denoted by $ A(n). $ \\
A subalgebra $ H $ of $ L $ is a vector subspace of $ L $ provided that  is closed under Lie bracket. 
By considering $ L $ as  an $ n $-dimensional Lie algebra over the field $ \mathbb{F}_{q} $  and  $ \lbrace x_1, \dots, x_n \rbrace$  as a basis of   $ L, $ 
 we may display $ a $ as a  linear  combinations of $  x_1, \dots, x_n $ for every element $ a\in L. $
In other words, we have  $ a=\alpha_1 x_1+\dots +\alpha_n x_n, $  where $ \alpha_i \in \mathbb{F}_q $ for all $ 1\leq i\leq n. $ 
Hence, there is $ q $ choices for every scaler and so the number of elements of $ L $ is equal to $ q^{n} $
and  the number of elements of an $ m $-dimensional subalgebra $ H $ is equal to $ q^{m} $ such that $ m\leq n. $
  An  ideal of $ L $ is a vector subspace of  $ L $ provided that  $ [I, L]\subseteq I. $  Every ideal is always a subalgebra, but a subalgebra 
  is not necessarily an ideal. For instance, let $ S $ be a subset of $ L. $ Then the centralizer of $ S $ in $ L $ is defined as  $ C_{L}(S)=\lbrace x\in L \mid [x, s]=0 ~~ \text{for all}~~ s\in S \rbrace $ which   is a subalgebra  and  not an ideal in general.
  There are  several examples for an ideal of $ L. $ One important example of an ideal of $ L $ is  the center of $ L $  defined as   $ Z(L)=\lbrace z\in L\mid [z, x]=0 \quad \text{for all} \quad  x\in L \rbrace.$ Note that $ L $ is abelian if and only if $Z(L)=L.$
Assume that 
$ I $ and $ J $ are two ideals of $ L. $ Then one can see $ [I, J]=\langle [x, y] \mid x\in I, y\in J \rangle $ is an ideal, which is named 
a product of ideals. By using it, we can construct another important example of an ideal of $ L. $ The derived subalgebra of $ L$ can be defined by  putting
   $ I=J=L. $ 
It is common to be denoted by $ L^{2} $ instead of $ [L, L]. $ One of known  classes of Lie algebra is Heisenberg Lie algebras. 
From \cite{Sti}, a Lie algebra $ L $ is called Heisenberg if $ L^{2}=Z(L) $ and $ \dim L^{2}=1 $, it is proved such Lie algebras are  odd-dimensional and have the following presentation
\begin{equation*}
L=H(m)=\langle  x_i, y_i, z \mid [ x_i, y_i]=z,~~ 1\leq i \leq m \rangle.
\end{equation*}  
Suppose that $ L_1 $ and $ L_2 $ are two Lie algebras. The (external) direct sum of them is the vector space 
$ L_1\oplus L_2=\lbrace (x_1, x_2) \mid x_1\in L_1, x_2\in L_2 \rbrace $  with the following operation 
\begin{equation*}
[(x_1, x_2), (y_1, y_2)]=([x_1, y_1], [x_2, y_2]). 
\end{equation*}
Put $ A=\lbrace (x_1, 0) \mid x_1 \in L_1 \rbrace$ and $ B=\lbrace (0, x_2) \mid x_2 \in L_2 \rbrace.$ It is easy to see that they are ideals of direct sum such that are isomorphic to $ L_1 $ and $ L_2, $ respectively. Also,  $ A\cap B=0 $ and $ [A, B]=0. $ 
Conversely, if $ I_1 $ and $ I_2 $ are two ideals of $ L $ such  that $ L=I_1+I_2 $ and $ I_1\cap I_2=0, $ then the map $ (x_1, x_2)\mapsto x_1+x_2 $ is an isomorphism of $ I_1\oplus I_2 $ to $ L. $ So, $ L $ is internal direct sum of $ I_1 $ and $ I_2. $ Note that 
$ [I_1, I_2] $ is contained in $ I_1\cap I_2 $ and so  $ [I_1, I_2]=0. $
 \\The lower central series of a Lie algebra $ L $ is a series with terms 
 \begin{equation*}
 L_0=L ~\text{and} ~ L_k=[L, L_{k-1}]~\text{for ~ all} ~k\geq 1.
 \end{equation*}
 Then the series $ L_0 \supseteq L_1 \supseteq L_2\supseteq \dots $ is decreasing.
 The Lie algebra  $ L $ is called nilpotent if $ L_m=0 $ for some $ m\geq 1. $   For a nilpotent Lie algebra, the smallest $ n $ such that $  L$ has a central series of length $ n $ is called the nilpotency class of $ L. $
 If the lower central series does not stop, then $ L $
 is called non-nilpotent Lie algebra.
 There  are various non-nilpotent
 Lie algebras. One of the non-nilpotent Lie algebras  is $ \langle x, y \mid [x, y]=x \rangle $ which plays an essential role   in the next. \\
 All notions and terminologies about Lie algebras are standard here and can be found in \cite{G, W}.

\section{ the Commutativity degree}
Let us start this section by definition of the commutativity degree of finite-dimensional Lie algebra $ L $ over the field $ \mathbb{F}_{q} $ as the following. 
\begin{defn}\label{def1.1}
Let $ L $ be a finite-dimensional  Lie algebra over the field $ \mathbb{F}_{q}. $ The commutativity degree of $ L $ is defined  as the ratio
\begin{equation*}
d(L)=\dfrac{| \lbrace (x, y)\in L\times L \mid [x, y]=0\rbrace |}{| L|^{2}}.
\end{equation*}
\end{defn}
It is clear that $ d(L)=1 $ if and only if $ L $ is abelian. So, we work on non-abelian Lie algebras. From now, all  Lie algebras are considered
 finite-dimensional non-abelian Lie algebra over the field $ \mathbb{F}_{q}. $\\
Put $ B(L)= \lbrace (x, y)\in L\times L \mid [x, y]=0\rbrace. $ We observe that $  | B(L)| =\sum_{x\in L} | C_{L}(x)|,  $ where 
$ C_{L}(x)=\lbrace y\in L \mid [x, y]=0 \rbrace. $ So, 
we can easily see that  
$ d(L)=\frac{1}{| L|^{2}}\sum_{x\in L} | C_{L}(x)|.$\\
For a Lie algebra $ L $ and every element $ x\in L, $ we denote $ ad_x $ as the map from $ L $ to $ L $
by the rule $ ad_x(y)=[x, y] $ for all $ y\in L. $
It is obvious that $ ad_x $ is a linear map.  
The following Lemma shows that   $ d(L) $ can be rewritten in term of   $ \operatorname{Im}ad_{x}$'s. 
\begin{lem}\label{lem2.2}
Let $ L $ be a Lie algebra. Then
\begin{equation*}
d(L)=\frac{1}{| L|}\sum_{x\in L} \frac{1}{|  \operatorname{Im}ad_{x}|}.
\end{equation*}
\end{lem}
\begin{proof}
Consider $ ad_{x} : L\rightarrow L. $ 
We have
$\dim L=\dim \operatorname{Im}ad_{x}+ \dim \ker ad_{x}. $ Since $ \ker ad_{x}=C_{L}(x), $ so $\dim L=\dim \operatorname{Im}ad_{x}+ \dim C_{L}(x). $ Hence 
\begin{align}\label{eq1.1}
| L |&=q^{\dim L}=q^{\dim \operatorname{Im}ad_{x}+ \dim C_{L}(x)}=|  \operatorname{Im}ad_{x}| | C_{L}(x)|. 
\end{align}
Therefore, the result follows by using  $ d(L)=\frac{1}{| L|^{2}}\sum_{x\in L} | C_{L}(x)|$ and  \eqref{def1.1}.
\end{proof}
By Lemma \ref{lem2.2}, we may compute the commutativity degree of  Heisenberg Lie algebras.
\begin{exa}\label{ex1}
Let  $ m $ be a positive integer and $ L=H(m)=\langle  x_i, y_i, z \mid [ x_i, y_i]=z,~~ 1\leq i \leq m \rangle. $
 Since $ \operatorname{Im}ad_x \subseteq L^{2} $ and $ \dim L^{2}=1, $ we have $ \operatorname{Im}ad_x = L^{2} $ for all $ x\in L\setminus Z(L). $ Hence 
$ |  \operatorname{Im}ad_{x}|=q $ for all $ x\in L \setminus Z(L) $ and 
\begin{align*}
d(L)&= \frac{1}{| L|}\sum_{x\in L} \frac{1}{|  \operatorname{Im}ad_{x}|}= \dfrac{|  Z(L) |  }{|  L |  }+\dfrac{1}{|  L |} \sum_{x\in L\setminus Z(L)} \dfrac{1}{|  \operatorname{Im}ad_{x}|}\cr
&= \dfrac{|  Z(L) |  }{|  L |  }+\dfrac{1}{|  L |}(|  L |- |  Z(L) |) \dfrac{1}{q}=\dfrac{q^{2m}+q-1}{q^{2m+1}}.
\end{align*} 
\end{exa}
In the following example, we compute the commutativity degree of the unique non-nilpotent Lie algebra of dimension $2$. 
\begin{exa}\label{ex2}
Let $ L=\langle  x_1, y_1 \mid [x_1, y_1]=x_1 \rangle. $
 It is clear that $ Z(L)=0 $ and so $ |  Z(L) |=1.  $ If $ x\neq 0, $ then $ \operatorname{Im}ad_x = L^{2}, $ which implies   $ |  \operatorname{Im}ad_{x}|=q. $ Thus
\begin{align*}
d(L)&= \frac{1}{| L|}\sum_{x\in L} \frac{1}{|  \operatorname{Im}ad_{x}|}= \dfrac{1 }{|  L |  }+\dfrac{1}{|  L |} \sum_{x\neq 0} \dfrac{1}{|  \operatorname{Im}ad_{x}|}= \dfrac{1 }{|  L |  }+\dfrac{1}{|  L |}(|  L |- |  Z(L) |)\frac{1}{q} \cr
&=\dfrac{q^{2}+q-1}{q^{3}}.
\end{align*} 
\end{exa}
   The commutativity degree of  a direct sum of two Lie algebras is given in the following.
\begin{prop}\label{pro3.1}
Let $ H $ and $ K $ be Lie algebras and $ L=H \oplus K. $ Then
\begin{equation*}
d(L)=d(H) d(K).
\end{equation*}
\end{prop}
\begin{proof}
Since
there is a natural bijection between $ B(L) $ and $  B(H)\times B(K), $ we have $ |B(L)|=|B(H)||B(K)|. $ So, 
\begin{align*}
| L |^{2}&=q^{\dim H+\dim K} q^{\dim H+\dim K}
= | H |^{2} | K |^{2}.
\end{align*}
Now, the result follows.
\end{proof}
A Lie algebra $ L $ is called stem if $ Z(L)\subseteq L^{2}. $ From \cite[Proposition 3.1]{J}, we  know that every finite-dimensional Lie algebra $ L $ is isomorphic to $ T\oplus A, $ in which
$ A $ is  abelian  and $ T $ is  stem. Now,   by Proposition \ref{pro3.1} we have $ d(L)=d(T) $. Therefore,  the commutativity degree of stem Lie algebras is important.
Thus, without loss of generality we may always assume that all Lie algebras are stem.
\\
The concept of isoclinism  is a useful tool in computation of the commutativity degree of groups. As we know two groups 
 $ G_1 $ and $ G_2 $ are called isoclinic provided that there are isomorphisms $ \alpha : G_1/Z(G_1)\rightarrow G_2/Z(G_2) $ and 
$ \beta : G_1^{\prime}\rightarrow G_2^{\prime} $ such that the following diagram commutes:
\begin{displaymath}
\xymatrix{ G_1/Z(G_1)\times  G_1/Z(G_1)  \ar[d]^{\alpha \times \alpha} \ar[r] &  G_1^{\prime}
 \ar[d]^{~~\beta}\\
G_2/Z(G_2)\times G_2/Z(G_2) \ar[r] & G_2^{\prime}}
\end{displaymath}
The horizontal  maps are defined by $ (\overline{x}, \overline{y})\mapsto [x, y]. $ This notion is denoted by $ G_1 \sim G_2. $ 
  Lescot \cite{Les2} proved that two isoclinic groups have the same commutativity degree. Now, it is interesting to see  what will happen for  isoclinic Lie algebras. We know that the definition of isoclinism  in 
 Lie algebras  the same as the group case 
by replacing $ G'_1 $ and $G'_2$ by $L^{2}_1$ and  $L^{2}_2.$ In the following theorem, we show that  two isoclinic Lie algebras $ L $ and $ M $ have the same commutativity degree similar to groups.
\begin{thm}\label{iso} 
Let $ L $ and $ M $ be the isoclinic Lie algebras.  Then $ d(L)=d(M). $
\end{thm}
\begin{proof}
By using \cite[Proposition 3.1]{J}, we have $ L=T_1\oplus A_1 $ and $M=T_2\oplus A_2 $ in which $ T_1 $ and $ T_2 $ are stem  and $ A_1 $ and $ A_2 $ are abelian. It is proved that \cite[Lemma 1]{Mon}
 $ L\sim T_1 $ and $ M\sim T_2 $ and hence   $T_1\sim T_2. $
Now, by using \cite[Theorem 11]{Mon} $ T_1\cong T_2$ and so $ d(T_1)=d(T_2). $ 
The result  follows by using Proposition \ref{pro3.1}.
\end{proof}
The concept of isoclinism   enables us to compute the commutativity degree of Lie algebras  more easily. The next example shows how it works.
\begin{exa}
Consider the following Lie algebras 
\[L=H(1)\oplus H(1)=\langle a_1, \dots, a_6\mid [a_1, a_2]=a_5, [a_3, a_4]=a_6 \rangle~~ \text{and}
\]
\[
M=L^{(2)}_{6,7}(0)=\langle x_1, \dots, x_6\mid [x_1, x_2]=x_5, [x_3, x_4]=x_5+x_6, [x_1, x_3]=x_6 \rangle
\]
over the field $ \mathbb{F}_{2}. $ We claim that $ d(L)=d(M)=\frac{25}{64}. $ By using the identity homomorphism $ \beta : L^{2}\longrightarrow M^{2} $ such that $ a_5\longmapsto x_5 $ and $ a_6\longmapsto x_5+x_6 $
and the 
 isomorphism $ \alpha : L/Z(L)\longrightarrow M/Z(M)$ given by
\begin{align*}
& \overline{a}_1\longmapsto \overline{x}_1+\overline{x}_2+\overline{x}_4, \quad \overline{a}_2\longmapsto \overline{x}_2, 
\quad \overline{a}_3\longmapsto \overline{x}_2+\overline{x}_3+\overline{x}_4, \quad \overline{a}_4\longmapsto \overline{x}_4,
\end{align*} 
 we can see  $L \sim M.$ Now, $d(L)=d(M)$ by Theorem \ref{iso}. On the other hand,
 $ d(L)= d(H(1))  d(H(1))= \frac{25}{64} $ by Example \ref{ex1} and Proposition \ref{pro3.1}. So, the result holds.
\end{exa}
It is known that  two groups with the same commutativity degree need  not   to be isoclinic (consider two groups $D_8$ and $Q_8$). Now, we have the same situation  for Lie algebras. 
\begin{exa}
Take the following  Lie algebras
\[ L_{5,5}=\langle x_1, \dots, x_5 \mid [x_{1},x_{2}] =x_{3}, [x_{1},x_{3}] =x_{5}, [x_{2},x_{4}] =x_{5} \rangle ~~\text{and} \]
\[ L_{5,7}=\langle x_1, \dots, x_5 \mid [x_{1},x_{2}] =x_{3}, [x_{1},x_{3}] =x_{4}, [x_{1},x_{4}]=x_{5} \rangle.\] 
First, we show that $ d(L_{5,5})= d(L_{5,7})=\frac{q^{3}+q^{2}-1}{q^{5}}. $ 
We know that $ \operatorname{Im}ad_{a}\subseteq L^{2} $ for every $ a\in L_{5,5} $ and $ \dim L^{2}_{5,5}=2 $ which implies
$ \dim \operatorname{Im}ad_{a}\leq 2.$ Also, we can write $ a=\alpha_1 x_1+\alpha_2 x_2+ \alpha_3 x_3+\alpha_4 x_4+\alpha_5 x_5+z$ for some $ \alpha_i \in \mathbb{F}_q $ such that   $ 1\leq i \leq 5 $ and $ z\in Z(L). $ Then $ \dim \operatorname{Im}ad_{a}= 2$ when $ \alpha_1 \neq 0. $
By counting the coefficients of $ a, $
the number of such elements is $ (q-1)q^{4} $ when $ \alpha_1\neq 0. $ If   $ \alpha_1 = 0 $ and $ \alpha_2 \neq 0, $ then
one can see $ \dim \operatorname{Im}ad_{a}= 2$ by a similar method and the number of them is  $ (q-1)q^{3}. $
Let $ \alpha_1= \alpha_2=0 $ and $ \alpha_3\neq 0. $ Then   $ \dim \operatorname{Im}ad_{a}= 1$ and their number is equal to 
$ (q-1)q^{2}. $ Finally, $ \dim \operatorname{Im}ad_{a}= 1$ and the number of these elements is $ (q-1)q $ when 
$ \alpha_1= \alpha_2=\alpha_3=0 $ and $ \alpha_4\neq 0. $ On the other hand, $ |Z(L_{5,5}) |=q $ thus by Lemma \ref{lem2.2}
\begin{align*}
&d(L_{5,5})=\frac{1}{| L|}\sum_{x\in L} \frac{1}{|  \operatorname{Im}ad_{x}|}=\frac{| Z(L)|}{| L|}+\frac{1}{| L|}\sum_{x\in L\setminus Z(L)} \frac{1}{|  \operatorname{Im}ad_{x}|}\cr
&=\frac{q}{q^{5}}+\frac{1}{q^{5}}(\frac{q^{3}-q}{q}+\frac{q^{5}-q^{3}}{q^{2}})
=\frac{q^{3}+q^{2}-1}{q^{5}}.
\end{align*}
Similarly, $ d(L_{5,7})=\frac{q^{3}+q^{2}-1}{q^{5}} $ is obtained.  
Since  $ |L^ {2}_{5,5}|\neq  |L^{2}_{5,7}|,$  they are not isoclinic.
\end{exa}

\section{Bounds on commutativity degree}
In this section, we give an upper bound and  a lower bound for the commutativity degree of a Lie algebra in term of $ q. $
As a result, we have  the upper bound $ \frac{5}{8}$ for the Lie algebra $ L $ (similar to the upper bound for the commutativity degree in groups).  It is worth emphasizing that there is no Lie algebra with the commutativity degree in the interval $ (\frac{5}{8}, 1). $
\\
We start with the upper bound for centerless Lie algebras, but before that we recall the breath of a Lie algebra used
 in the next theorems. \\
For  any Lie algebra $ L, $  the breath of $ L $ is denoted by $ b(L) $  defined as the set $ b(L)=\max \lbrace \dim \operatorname{Im}ad_x \mid x\in L \rbrace.$ This concept plays an essential role in  classification of Lie algebras and other  issues (see \cite{B} for more details).
\begin{thm}\label{Bon1}
Let $ L $ be an $ n $-dimensional Lie algebra with $ Z(L)=0. $ Then 
\begin{equation*}
 d(L)\leq \dfrac{q^{n}+q-1}{q^{n+1}}.
\end{equation*}
The equality holds if and only if $ \dim L^{2}=1. $
\end{thm}
\begin{proof}
  If $ x\neq 0, $ then $ \operatorname{Im}ad_x  $ is a non-zero and 
 $\operatorname{Im}ad_x\subseteq  L^{2} $ and so $ \dim \operatorname{Im}ad_x \geq 1.$ Thus 
   $ \frac{1}{|\operatorname{Im}ad_x|}\leq \frac{1}{q} $ for all non-zero elements $ x. $ Hence by  Lemma \ref{lem2.2}
\begin{align*}\label{eq3.3}
d(L)&=\frac{1}{| L|}\sum_{x\in L} \frac{1}{|  \operatorname{Im}ad_{x}|}=\dfrac{1}{| L|}+\dfrac{1}{| L|} \sum_{x\neq 0}\dfrac{1}{|\operatorname{Im}ad_x|} \cr
&\leq \dfrac{1}{q^{n}}+\dfrac{1}{q^{n}} \dfrac{q^{n}-1}{q}
= \dfrac{q^{n}+q-1}{q^{n+1}}.
\end{align*}
Assume that $ \dim L^{2}=1. $ Then $ b(L)=1 $ by \cite[Theorem 3.2.1]{B} and so $ |\operatorname{Im}ad_x|=q $ for all 
$ x\in L\setminus Z(L). $ Therefore, $ d(L)= \dfrac{q^{n}+q-1}{q^{n+1}}$ by Lemma \ref{lem2.2}. Conversely, let the equality hold. We claim that $ \dim L^{2}=1. $ By contrary, suppose that  $ \dim L^{2}\neq 1. $ Since $ L $ is non-abelian, then $ \dim L^{2}\geq 2. $
So, there exists an element $ x $ such that $ |\operatorname{Im}ad_x|> q$ by \cite[Theorem 3.2.1]{B}. Hence 
$ d(L)< \dfrac{q^{n}+q-1}{q^{n+1}},$ which is a contradiction.
\end{proof}
The case that $ Z(L)\neq 0 $ is stated as the following.
\begin{thm}\label{Bon11}
Let $ L $ be a Lie algebra with $ \dim Z(L)\geq 1. $  Then 
\begin{equation*}
 d(L)\leq \dfrac{q^{2}+q-1}{q^{3}}.
\end{equation*}
\end{thm}
\begin{proof}
By  Lemma \ref{lem2.2}, we have
\begin{align*}
d(L)&=\frac{1}{| L|}\sum_{x\in L} \frac{1}{|  \operatorname{Im}ad_{x}|}
=\dfrac{1}{| L|}\big(| Z(L)|+\sum_{x\in L\setminus Z(L)}\dfrac{1}{|\operatorname{Im}ad_x|}\big)\cr
&\leq \dfrac{|Z(L)|}{|L|}+\dfrac{1}{q}-\dfrac{|Z(L)|}{|L| q}=\dfrac{1}{q}+\dfrac{|Z(L)|}{|L|}(1-\dfrac{1}{q}).
\end{align*}
One the other hand, we know that $ \dim L/Z(L)\geq 2 $ for all non-abelian Lie algebras.
  Thus 
\begin{equation*}\label{eq3.2}
d(L)\leq \dfrac{q^{2}+q-1}{q^{3}}.
\end{equation*}
as required.
\end{proof}
\begin{cor}\label{corBB}
Let $ L $ be a  Lie algebra.   Then 
\begin{equation}\label{BB}
d(L)\leq \dfrac{q^{2}+q-1}{q^{3}}.\end{equation}
\end{cor}
\begin{proof}
Assume that  $ \dim L=n. $ If $ Z(L)\geq 1, $ then $ d(L)\leq \frac{q^{2}+q-1}{q^{3}}$ by Theorem \ref{Bon11}.
If $ Z(L)=0, $ then $ d(L)\leq  \frac{q^{n}+q-1}{q^{n+1}}$ by Theorem \ref{Bon1}. It is enough to show that 
$  \frac{q^{n}+q-1}{q^{n+1}} \leq \frac{q^{2}+q-1}{q^{3}}. $ 
 We know $ q^{3}(q-1)\leq q^{n+1} (q-1) $ for all $ n\geq 2, $ which implies 
$ \frac{q-1}{q^{n+1}}\leq \frac{q-1}{q^{3}}.$ Hence $  \frac{q^{n}+q-1}{q^{n+1}} \leq \frac{q^{2}+q-1}{q^{3}}. $ 
\end{proof}
 By  the  breath of a Lie algebra,  we give a lower bound
 for the commutativity degree of Lie algebras. Let us state the following Lemma from \cite{B} which can be used in the proof of Theorem \ref{low}.
\begin{lem}\cite[Theorem 3.1.9]{B}\label{b}
Let $ L $ be a Lie algebra such that $ b(L)=n$, ($n>0$). Then $ \dim(L/Z(L)) \geq n+1.  $
\end{lem}
\begin{thm}\label{low}
Let $ L $ be a Lie algebra such that  $ \dim L/Z(L)=t. $ Then 
\begin{equation*}
d(L)\geq \frac{1}{q^{t}}+\frac{1}{q^{t-1}}-\frac{1}{q^{2t-1}}.
\end{equation*}
\end{thm}
\begin{proof}
 Let $ \dim L/Z(L)=t. $ Theorem \ref{b} shows that $ \dim L/Z(L)\geq \dim \operatorname{Im}ad_x+1 $ for all $ x\in L, $ so we have $ \frac{1}{|\operatorname{Im}ad_x|}\geq\frac{1}{q^{t-1}}.$ Hence 
  by  Lemma \ref{lem2.2} 
\begin{align*}\label{eqq3.4}
d(L)&=\frac{1}{| L|}\sum_{x\in L} \frac{1}{|  \operatorname{Im}ad_{x}|} =\dfrac{1}{| L|}\big( | Z(L)|+\sum_{x\in L\setminus Z(L)}\dfrac{1}{|\operatorname{Im}ad_x|}\big)\cr
&\geq \frac{| Z(L)|}{| L|}+\dfrac{( | L|- | Z(L)|) }{|L|}\dfrac{1}{q^{t-1}}= \dfrac{1}{q^{t}}+\dfrac{1}{q^{t-1}}-\dfrac{1}{q^{2t-1}}.
\end{align*}
\end{proof}
The next corollary shows all central quotients of a Lie algebra $ L $ that $ d(L) $ attains the bound \eqref{BB} has exactly $ 2 $-dimensional.
\begin{cor}\label{cor3.4}
Let $ L $ be a   Lie algebra. Then  $ \dim L/Z(L)=2$ if and only if $ d(L) = \frac{q^{2}+q-1}{q^{3}}. $
\end{cor}
\begin{proof}
Since $ \dim L/Z(L)=2,$  then $ \frac{q^{2}+q-1}{q^{3}} \leq d(L) $ by  
Theorem \ref{low}. On other hand, $ d(L)\leq \frac{q^{2}+q-1}{q^{3}} $ by Corollary
\ref{BB}. 
  So, $ d(L) = \frac{q^{2}+q-1}{q^{3}}. $ 
Conversely, let $ x\notin Z(L). $ Then $ | L |/ | C_{L}(x) |\geq q.$ Hence
\begin{align*}
\frac{q^{2}+q-1}{q^{3}}=d(L)=\frac{1}{| L|^{2}}\sum_{x\in L} | C_{L}(x)|&\leq \frac{1}{| L|^{2}} (| Z(L) || L |+(| L |-| Z(L) |)\dfrac{| L |}{q} )\cr
&=\dfrac{1}{q}+\dfrac{| Z(L) |}{| L |}(1-\dfrac{1}{q}).
\end{align*}
Thus $ \frac{| L |}{| Z(L) |}\leq q^{2}. $ Since $ L $ is non-abelian, we have $ \dim L/Z(L)=2.$
\end{proof}
In the next theorem, we obtain the upper bound $ \frac{5}{8}$ and show that it is just valid for the commutativity degree of all Lie algebras over the field $ \mathbb{F}_{2}. $ Hence there is no Lie algebra $L$ with $d(L)>\frac{5}{8}.$
\begin{thm}
Let $ L $ be an $ n $-dimensional  Lie algebra. Then $  d(L)\leq \frac{5}{8}. $ Moreover, $ d(L)=\frac{5}{8} $ if and only if $ L $ is isomorphic to one of the following Lie algebras over the field $ \mathbb{F}_2 $
\[  \langle x,y, z \mid [x, y]=z \rangle = H(1)\oplus A(n-3) \]
or 
\[ \langle x,y \mid [x, y]=x \rangle \oplus A(n-2).
\]
\end{thm}
\begin{proof}
By  Corollary \ref{corBB}, we have $ d(L)\leq \frac{q^{2}+q-1}{q^{3}}. $ Put $ f(q)= \frac{q^{2}+q-1}{q^{3}}.$
Since $ f(q) $ is a decreasing  function for all $ q\geq 2, $ we have $  d(L)\leq \frac{5}{8}. $\\
If $ L $ is isomorphic to one of the above  Lie algebras, 
then $ d(L)=\frac{5}{8} $  by  Proposition \ref{pro3.1}, Examples \ref{ex1} and \ref{ex2} for $ q=2. $ 
Conversely, let $ d(L)=\frac{5}{8}. $  Note that $ d(L)=\frac{5}{8} $ when $ q=2 $ by the first proof.
Since $ L=T\oplus A $ where $ T $ is  stem  and $ A $ is abelian  by   \cite[Proposition 3.1]{J}, it is enough to determine the structure of $ T. $  Also, $\dim T/Z(T)=2$ by  Corollary \ref{cor3.4} thus 
$ \dim T^{2}\leq 1 $ by \cite[Lemma 14]{Mon}. On the other hand, $ T $ is non-abelian, so $ \dim T^{2}=1. $
 Since $ Z(T)\subseteq T^{2} $ and $ \dim T^{2}=1, $ we have $ \dim Z(T)=0 $ or $ \dim Z(T)=1. $ 
Therefore $ \dim T=2 $ or $ 3 $  and so 
 $ T $ is isomorphic to $ \langle x,y, z \mid [x, y]=z \rangle \cong H(1)$ or $ \langle x,y \mid [x, y]=x \rangle  $ by  \cite[Theorem 3.1 and section 3.2.1]{W},  Examples \ref{ex1} and \ref{ex2}.  Hence $ L $ is isomorphic to one of the Lie algebras over the field $ \mathbb{F}_2 $
\[  \langle x,y, z \mid [x, y]=z \rangle \cong H(1)\oplus A(n-3), \]
or 
\[ \langle x,y \mid [x, y]=x \rangle \oplus A(n-2).
\]
\end{proof}
Now, we are in a position to compare $ d(L) $ and $ d(H),$ where $H$ is a subalgebra of $L$.
\begin{lem}\label{lem4.2}
Let $ L $ be a  Lie algebra and $ H $ be a subalgebra of $ L. $ Then
$  \frac{| H|}{| C_H(x)|}\leq \frac{| L|}{| C_L(x)|}$
for all $ x\in L. $
\end{lem}
\begin{proof}
It is clear that $H+C_L(x)\subseteq L$ for all $ x\in L. $ 
Hence $|  H+C_L(x)| \leq | L|   $ and so $ |  H| |C_L(x)| / |  H \cap C_L(x)| \leq | L|.   $ Since 
$  H \cap C_L(x)= C_H(x),$ we have $ |  H| |C_L(x)| / |   C_H(x)| \leq | L| $ as required.
\end{proof}
The following proposition gives a comparison between   $ d(H) $  and $ d(L) $ for every subalgebra  $ H $ of $ L. $ 
\begin{prop}\label{pro2.5}
Let $ L $ be a Lie algebra and $ H $ be a subalgebra of $ L. $ Then  $ \frac{| H |^{2}}{| L |^{2}} d(H)\leq d(L) \leq d(H). $ If  $L=H+Z(L), $ then $ d(L)=d(H). $
\end{prop}
\begin{proof}
 By Lemma \ref{lem4.2}
\begin{equation}\label{eqsb1}
\sum_{x\in L} | C_{L}(x) |\leq \dfrac{| L |}{| H |}\sum_{x\in L} | C_{H}(x) |.
\end{equation}
If $ x\in L $ and $ y\in  C_{H}(x),$ then $ [x, y]=0. $ Thus $ x \in  C_{L}(y) $ and we can rewrite \eqref{eqsb1} as follows:
\begin{equation}\label{eqsb2}
\sum_{x\in L} | C_{L}(x) |\leq \dfrac{| L |}{| H |}\sum_{y\in H} | C_{L}(y) |.
\end{equation}
 By Lemma \ref{lem4.2} and \eqref{eqsb2} we have
\begin{equation}\label{eqsb22}
\sum_{x\in L} | C_{L}(x) |\leq \dfrac{| L |^{2}}{| H |^{2}} \sum_{y\in H} | C_{H}(y) |.
\end{equation}
Consequently,
\begin{align*}
d(L)=\dfrac{1}{| L |^{2}}\sum_{x\in L} | C_{L}(x) | \leq \dfrac{1}{| L |^{2}}\dfrac{| L |^{2}}{| H |^{2}} \sum_{y\in H} | C_{H}(y) |=\dfrac{1}{| H |^{2}} \sum_{y\in H} | C_{H}(y) |=d(H).
\end{align*}
Since $| C_{H}(x) | \leq | C_{L}(x) |, $ we have also
\begin{align*}
d(L)&=  \dfrac{1}{| L |^{2}}\sum_{x\in L} | C_{L}(x) | \geq  \dfrac{1}{| L |^{2}}\sum_{x\in L} | C_{H}(x) |  \geq \dfrac{1}{| L |^{2}}\sum_{x\in H} | C_{H}(x) | =\dfrac{| H |^{2} }{| H |^{2} | L |^{2}}  \sum_{x\in H} | C_{H}(x) |  \cr
&=\dfrac{| H |^{2}}{| L |^{2}} d(H) .
\end{align*}
Therefore $ \frac{| H |^{2}}{| L |^{2}} d(H)\leq d(L) \leq d(L).  $\\
Let $ L=H+Z(L). $ Then  $ Z(L)=\cap_{x\in L} C_L(x) $ and $ H+ C_L(x)\subseteq L $  imply that  $ L=H+ C_L(x). $ Hence 
$  \frac{| H|}{| C_H(x)|}=\frac{| L|}{| C_L(x)|}$ for all $ x\in L $ by  the proof of Lemma \ref{lem4.2}.
Therefore  the inequality change into equality in  \eqref{eqsb22} and so $ d(L)=d(H). $ 
\end{proof}
 Now, we compare $ d(L) $ with $ d(N) $ and $ d(L/N),$ where   $ N $ is  an ideal of $ L.$
\begin{prop}\label{pro2.6}
Let $ L $ be a  Lie algebra and $ N $ be an ideal of $ L. $ Then 
\begin{equation*}
 d(L) \leq d(\dfrac{L}{N}) d(N).
\end{equation*}
If $N \cap L^{2}=0, $ then the equality holds.
\end{prop}
\begin{proof}
Let $y \in L. $  By  the second isomorphism theorem and $ C_N(y)= C_L(y)\cap N, $  we have 
$ \frac{C_L(y)}{C_N(y)} \cong \frac{C_L(y)+ N}{N}\subseteq C_{L/N}(y+N).$ Hence 
\begin{align}\label{eq3.9}
\sum_{y\in L}|C_L(y) |&=\sum_{S\in L/N} \sum_{y\in S}\dfrac{|C_L(y) |}{|C_N(y) |}|C_N(y) |=\sum_{S\in L/N} \sum_{y\in S}\dfrac{|C_L(y) |}{|C_L(y)\cap N |}|C_N(y) |\cr
&=\sum_{S\in L/N} \sum_{y\in S}\dfrac{|C_L(y)+N |}{|N |}|C_N(y) |\leq  \sum_{S\in L/N} \sum_{y\in S}|C_{L/N}(y+N) ||C_N(y) |\cr
&=  \sum_{S\in L/N} |C_{L/N}(S) | \sum_{y\in S} |C_N(y) |=  \sum_{S\in L/N}|C_{L/N}(S) | \sum_{y\in S}|\lbrace x\in N \mid [x, y]=0 \rbrace |\cr
&=  \sum_{S\in L/N} |C_{L/N}(S) | \sum_{x\in N}|C_L(x) \cap S |.\cr
\end{align}
If $ C_L (x) \cap S=\emptyset, $ then $ |C_L(x) \cap S |<  |C_N(x) |.$
If $ C_L (x) \cap S\neq \emptyset, $ then  there exist $ x_0 \in  C_L (x) \cap S. $ So, we have  $ S=x_0+N $ 
\begin{align*}
C_L (x) \cap S&= C_L (x) \cap (x_0+N)=(x_0+C_L (x)) \cap (x_0+N)=(C_L (x) \cap N)+x_0\cr
&=C_N(x)+x_0.
\end{align*}
Hence $ |C_L(x) \cap S |=|C_N(x)+x_0 |=|C_N(x) |.$ So, in any case
\begin{equation}\label{eq3.10}
|C_L(x) \cap S |\leq |C_N(x) |.
\end{equation}
Finally, by using  \eqref{eq3.9} and  \eqref{eq3.10}
\begin{align*}
d(L) &=\dfrac{1}{|L|^{2}}  \sum_{S\in L/N} |C_{L/N}(S) | \sum_{x\in N}|C_L(x) \cap S |\cr
&\leq \dfrac{|N|^{2}}{|L|^{2}}\sum_{S\in L/N} |C_{L/N}(S) | \dfrac{1}{|N|^{2}}\sum_{x\in N} |C_N(x) |= d(L/N) d(N).
\end{align*}
Assume that $ x+N \in C_{L/N}(y+N). $ Thus $ [x, y]\in N. $ If $ N \cap  L^{2}=0,$ then $ [x, y]=0 $ and so $ x+N \in \frac{C_L(y)+ N}{N}. $ 
Hence $  C_{L/N}(y+N)\subseteq \frac{C_L(y)+ N}{N}.$
Since $  \frac{C_L(y)+ N}{N}\subseteq C_{L/N}(y+N),$  we have   $\frac{C_L(y)+ N}{N}= C_{L/N}(y+N).$ We can also see that $ C_L (x) \cap S\neq \emptyset $ for all
$ x\in N $ and for $ S\in \frac{L}{N}. $ Therefore  all inequalities change into equalities. 
\end{proof}
\begin{prop}
Let $ L $ be a  Lie algebra   and  $ M, $ $ N $  be ideals of $ L $ such that $ N\subseteq M. $ Then 
$ d(\frac{L}{N})\leq d(\frac{L}{M}). $
\end{prop}
\begin{proof}
Since $d(\frac{L}{M})=d(L/N/M/N) $ by  the third isomorphism theorem, we have 
\begin{equation*}
d(\frac{L}{M})=d(\dfrac{L/N}{M/N}) \geq \dfrac{d(L/N)}{d(M/N)}\geq d(\frac{L}{N})
\end{equation*}
by Proposition \ref{pro2.6} and the fact  that $ 0< d(M/N) \leq 1.$ 
\end{proof}

\end{document}